\documentclass [12pt] {article}
\usepackage{latexsym,amsmath,amsfonts,amssymb,amsthm}
\newtheorem{thm}{Theorem}[section]

\newtheorem{cor}{Corollary}[section]
\newtheorem{conj}{Conjecture}[section]
\makeatletter
\@addtoreset{equation}{section}
\date{}{}
\def\colon{{:}\;}
\def\Z{\Bbb Z}

\def\N{\Bbb N}

\def\C{\Bbb C}
\def\al{\alpha}
\def\l{\left}
\def\r{\right}
\def\bg{\bigg}
\def\({\bg(}
\def\){\bg)}
\def\[{\bg\lfloor}
\def\]{\bg\rfloor}
\def\t{\mbox}
\def\f{\frac}

\def\em{\emptyset}
\def\se {\subseteq}

\def\sm{\setminus}

\def\bi{\binom}
\def\eq{\equiv}
\def\cs{\ldots}
\def\ls{\leqslant}
\def\gs{\geqslant}
\def\mo{\t{\rm mod}}

\def\sign{\t{\rm sign}}

\def\da{\delta}

\def\qed{}
\def\Proof{\noindent{\it Proof}}

\def\Remark{\medskip\noindent{\it  Remark}}
\def\Ack{\medskip\noindent {\bf Acknowledgments}}
\begin{document}
\hbox{Finite Fields Appl. 14(2008), no.\,2, 470--481.}
\bigskip
\centerline {\bf On Value Sets of Polynomials over a Field}
\bigskip
\centerline {Zhi-Wei Sun}
\bigskip
\centerline{Department of Mathematics, Nanjing University}
\centerline{Nanjing 210093, People's Republic of China}
\centerline{zwsun\char'100nju.edu.cn} \centerline{\tt
http://math.nju.edu.cn/$\sim$zwsun}
\bigskip
\abstract {Let $F$ be any field. Let $p(F)$ be the characteristic of $F$
if $F$ is not of characteristic zero, and let $p(F)=+\infty$ otherwise.
Let $A_1,\ldots,A_n$ be finite nonempty subsets of $F$, and let
$$f(x_1,\ldots,x_n)=a_1x_1^k+\cdots+a_nx_n^k+g(x_1,\ldots,x_n)\in F[x_1,\ldots,x_n]$$
 with $k\in\{1,2,3,\ldots\}$, $a_1,\ldots,a_n\in F\sm\{0\}$ and $\deg g<k$.
 We show that
 $$|\{f(x_1,\ldots,x_n):\,x_1\in A_1,\ldots,x_n\in A_n\}|
 \gs\min\bg\{p(F),\,\sum_{i=1}^n\l\lfloor\f{|A_i|-1}k\r\rfloor+1\bg\}.$$
 When $k\gs n$ and $|A_i|\gs i$ for $i=1,\ldots,n$, we also have
 \begin{align*} &|\{f(x_1,\ldots,x_n):\,x_1\in A_1,\ldots,x_n\in A_n,\ \t{and}\ x_i\not=x_j\ \t{if}\ i\not=j\}|
 \\&\qquad\qquad\gs\min\bg\{p(F),\,\sum_{i=1}^n\l\lfloor\f{|A_i|-i}k\r\rfloor+1\bg\};
 \end{align*}
consequently, if $n\gs k$ then for any finite subset $A$ of $F$ we have
$$|\{f(x_1,\ldots,x_n):\, x_1,\ldots,x_n\in A,\,\t{and}\ x_i\not=x_j\ \t{if}\ i\not=j\}|
\gs\min\{p(F),\,|A|-n+1\}.$$
In the case $n>k$, we propose a further conjecture
which extends the Erd\H os-Heilbronn conjecture in a new direction.}

\footnote{
\newline{\indent \quad Supported by
the National Science Fund (Grant No. 10425103) for Distinguished Young Scholars in China.}}

\section{Introduction}\par

 For a field $F$ we use $p(F)$ to denote the additive order of the multiplicative
identity of $F$, which is either infinite or a prime.
A field $F$ is said to have characteristic zero if $p(F)=+\infty$,
and have characteristic $p$ if $p(F)$ is a prime $p$.

By the Chevalley--Warning theorem (cf. pp. 50--51 of \cite{N}),
for any polynomial $P(x_1,\ldots,x_n)$ over a finite field $F$, if $n>\deg P$ then
the characteristic of $F$ divides the number of solutions to the equation $P(x_1,\ldots,x_n)=0$
over $F^n$. However, this says nothing about the solvability of the equation over $F^n$
unless there is an obvious solution.

Given a field $F$, we consider polynomials of the form
\begin{equation}\label{eq-1.1}
f(x_1,\ldots,x_n)=a_1x_1^k+\cdots+a_nx_n^k+g(x_1,\ldots,x_n)\in F(x_1,\ldots,x_n),
\end{equation}
where
\begin{equation}\label{eq-1.2}
k\in\Z^+=\{1,2,3,\ldots\},\ a_1,\ldots,a_n\in F^*=F\sm\{0\}\  \t{and}\  \deg g<k.
\end{equation}
What can we say about the solvability of the equation $f(x_1,\ldots,x_n)=0$ over $F^n$?

Let $F$ be the field $F_p=\Z/p\Z$ with $p$ a prime, and assume (1.1) and (1.2).
In 1956 Carlitz \cite{C} proved that the equation $f(x_1,\ldots,x_n)=0$ has a solution in $F_p^n$
when $k\mid p-1$ and $n\gs k$.
In 2006 Felszeghy \cite{F} extended this result by showing that the equation is solvable
if $k<p$ and $n\gs\lceil(p-1)/\lfloor(p-1)/k\rfloor\rceil$.
(For a real number $\al$, $\lceil\alpha\rceil$ denotes the least integer not smaller than $\al$
while $\lfloor\al\rfloor$ represents the largest integer not exceeding $\al$.)
Note that the equation $f(x_1,\ldots,x_n)=0$ over $F_p^n$ is solvable if and only if the value set
$$\{f(x_1,\ldots,x_n):\,x_1,\ldots,x_n\in F_p\}$$
contains $0$.
In 1959 Chowla, Mann and Straus (cf. pp. 60-61 of \cite{N}) used Vosper's theorem (cf.
pp. 52-57 of \cite{N}) to deduce that if $p>3$, $1<k<(p-1)/2$ and $k\mid p-1$, then
$$|\{a_1x_1^k+\cdots+a_nx_n^k:\, x_1,\ldots,x_n\in F_p\}|\gs\min\l\{p,\,(2n-1)\f{p-1}k+1\r\}.$$

Let $F_q$ be the finite field of $q$ elements where $q>1$ is a
prime power. In 1993 Wan, Shiue and Chen \cite{WSC} showed that if
$P(x)$ is a polynomial over $F_q$ and $l\in\N=\{0,1,2,\ldots\}$ is
the least nonnegative integer with $\sum_{x\in F_q}P(x)^l\not=0$
then $|\{P(x):\,x\in F_q\}|\gs l+1$; in 2004 Das \cite{D} extended
this to multi-variable polynomials over $F_q$. By modifying the
proof of Theorem 1.5 of \cite{D} slightly, one gets the following
assertion: If $P(x_1,\ldots,x_n)\in F_q[x_1,\ldots,x_n]$,
$\em\not=S\se F_q^n$, and $l$ is the smallest nonnegative integer
with $\sum_{(x_1,\ldots,x_n)\in S}P(x_1,\ldots,x_n)^l\not=0$, then
we have $|\{P(x_1,\ldots,x_n):\,(x_1,\ldots,x_n)\in S\}|\gs l+1$.
Here the lower bound depends heavily on values of
$P(x_1,\ldots,x_n)$.

In this paper we investigate two kinds of value sets of a polynomial in the form (1.1).

Here is our first theorem which includes Felszeghy's result as a special case.

\begin{thm}\label{Theorem 1.1} Let $f(x_1,\ldots,x_n)$ be a polynomial over a field $F$ given by $(1.1)$ and $(1.2)$.
Then, for any finite nonempty subsets $A_1,\ldots,A_n$ of $F$, we have
\begin{equation}\label{eq-1.3}|\{f(x_1,\ldots,x_n):\ x_1\in A_1,\ldots,x_n\in A_n\}|
\gs\min\l\{p(F),\,\sum_{i=1}^n\l\lfloor\f{|A_i|-1}k\r\rfloor+1\r\}.
\end{equation}
\end{thm}

This result will be proved in Section 2
where we also give an example to show that
the inequality (1.3) is sharp when $F$ is an algebraically closed field.

\begin{cor}\label{Corollary 1.1} {\rm (Cauchy--Davenport Theorem)}
Let $A_1,\ldots,A_n$ be finite nonempty subsets of a field $F$. Then
$$|A_1+\cdots+A_n|\gs\min\{p(F),\,|A_1|+\cdots+|A_n|-n+1\},$$
where the sumset $A_1+\cdots+A_n$ is given by
$$A_1+\cdots+A_n=\{x_1+\cdots+x_n:\, x_1\in A_1,\ldots,x_n\in A_n\}.$$
\end{cor}
\Proof. Simply apply Theorem 1.1 with $f(x_1,\ldots,x_n)=x_1+\cdots+x_n$. \qed

\Remark\ 1.1. The original Cauchy--Davenport theorem (see p.\,44 of \cite{N} or p.\,200 of \cite{TV})
is Corollary 1.1 in the case $n=2$ and $F=\Z/p\Z$ with $p$ a prime.

\begin{cor}\label{Corollary 1.2} Let $F$ be a field of characteristic zero, and assume $(1.1)$ and $(1.2)$.
Then, for any finite nonempty subset $A$ of $F$, we have
$$|\{f(x_1,\ldots,x_n):\, x_1,\ldots,x_n\in A\}|
\gs n\l\lfloor\f{|A|-1}k\r\rfloor+1.$$
\end{cor}
\Proof. Just apply Theorem 1.1 with $A_1=\cdots=A_n=A$. \qed

\begin{cor}\label {Corollary 1.3} Let $F$ be a field with prime characteristic $p$, and let
$f(x_1,\ldots,x_n)$ be given by $(1.1)$ and $(1.2)$. If $A$ is a finite subset of $F$
satisfying $\lfloor(|A|-1)/k\rfloor\gs(p-1)/n$, then
$$|\{f(x_1,\ldots,x_n):\, x_1,\ldots,x_n\in A\}|\gs p.$$
\end{cor}
\Proof. It suffices to apply Theorem 1.1 with $A_1=\cdots=A_n=A$. \qed

\Remark\ 1.2. In the case $A=F=F_p$, Corollary 1.3 yields Felszeghy's result.
\medskip

Now we state our second theorem.

\begin{thm}\label{Theorem 1.2} Let $f(x_1,\ldots,x_n)$ be a polynomial over a field $F$ given by
$(1.1)$ and $(1.2)$ with $n\ls k=\deg f$. And let
$A_1,\ldots,A_n$ be finite subsets of $F$ with $|A_i|\gs i$
for $i=1,\ldots,n$. Then, for the restricted value set
\begin{equation}\label{eq-1.4}
V=\{f(x_1,\ldots,x_n):\ x_1\in A_1,\ldots,x_n\in A_n,\ \t{and}\ x_i\not=x_j\ \t{if}\ i\not=j\},
\end{equation}
we have
\begin{equation}\label{eq-1.5}
|V|\gs\min\l\{p(F),\,\sum_{i=1}^n\l\lfloor\f{|A_i|-i}k\r\rfloor+1\r\}.
\end{equation}
\end{thm}

\begin{cor}\label{Corollary 1.4} Let $A$ be a finite subset of a field $F$, and let
$f(x_1,\ldots,x_n)$ be a polynomial given by $(1.1)$ and $(1.2)$.
Write $|A|=kq+r$ with $q,r\in\N$ and $r<k$. Then we have
\begin{equation}\label{eq-1.6}\aligned&|\{f(x_1,\ldots,x_n):\,
x_1,\ldots,x_n\in A,\,\t{and}\ x_i\not=x_j\ \t{if}\ i\not=j\}|
\\&\ \ \gs
\begin{cases}\min\{p(F),\,n(q-1)+\min\{n,r\}+1\}&\t{if}\ n\ls k,
\\\min\{p(F),\,|A|-n+1\}&\t{if}\ n\gs k.
\end{cases}
\endaligned
\end{equation}
\end{cor}

In Section 3 we shall prove Theorem 1.2 and Corollary 1.4, and give an example to illustrate that
the inequality (1.5) is essentially best possible when $F$ is an algebraically closed field.

In 1964 Erd\H os and Heilbronn \cite{EH} conjectured that if $A$ is a subset of $\Z/p\Z$ with $p$ a prime
then
$$|\{x_1+x_2:\, x_1,x_2\in A\ \t{and}\ x_1\not=x_2\}|\gs\min\{p,\,2|A|-3\}.$$
Thirty years later this deep conjecture was confirmed by Dias da Silva and Hamidoune \cite{DH}
who used the representation theory of groups to show that if $A$ is a finite subset of a field $F$ then
\begin{align*}&|\{x_1+x_2+\cdots+x_n:\, x_1,\ldots,x_n\in A,\ \t{and}\ x_i\not=x_j\ \t{if}\ i\not=j\}|
\\&\qquad\qquad\quad\gs\min\{p(F),\,n(|A|-n)+1\}.
\end{align*}
This suggests that Corollary 1.4 in the case $n>k$ might be further improved. In Section 4 we will discuss
our following conjecture which extends the Dias da Silva--Hamidoune result in a new way.
(The reader may consult \cite{ANR}, \cite{HS}, \cite{K} and \cite{S}
 for other generalizations of the Erd\H os-Heilbronn conjecture.)

\begin{conj}\label{Conjecture 1.1} Let $f(x_1,\ldots,x_n)$
be a polynomial over a field $F$ given by $(1.1)$ and $(1.2)$, and let
$A$ be any finite subset of $F$. Provided $n>k$, we have
\begin{equation}\label{eq-1.7}\aligned&|\{f(x_1,\ldots,x_n):\,
 x_1,\ldots,x_n\in A,\ \t{and}\ x_i\not=x_j\ \t{if}\ i\not=j\}|
\\\ &\gs\min\l\{p(F)-\da,\ \f{n(|A|-n)}k-k\l\{\f nk\r\}\l\{\f{|A|-n}k\r\}+1\r\},
\endaligned
\end{equation}
where $\{\al\}$ denotes the fractional part $\al-\lfloor\al\rfloor$
of a real number $\al$, and
$$\da=
\begin{cases}1&\t{if}\ n=2\ \t{and}\ a_1=-a_2,\\0&\t{otherwise}.
\end{cases}$$
\end{conj}

\section{Proof of Theorem 1.1}\par

We need a useful tool of algebraic nature.
\medskip

\noindent{\bf Combinatorial Nullstellensatz} (Alon \cite{A}).
{\it Let $A_1,\cs,A_n$ be finite subsets of a field
$F$, and let $P(x_1,\cs,x_n)\in F[x_1,\cs,x_n]$.
Suppose that $\deg P=k_1+\cdots+k_n$ where $0\ls k_i<|A_i|$ for $i=1,\ldots,n$.
Then $P(x_1,\ldots,x_n)\not=0$ for some $x_1\in A_1,\ldots,x_n\in A_n$ if
$$[x_1^{k_1}\cdots x_n^{k_n}]P(x_1,\ldots,x_n)\not=0,$$
where $[x_1^{k_1}\cdots x_n^{k_n}]P(x_1,\ldots,x_n)$ denotes the coefficient
of $x_1^{k_1}\cdots x_n^{k_n}$ in $P(x_1,\ldots,x_n)$.}

\medskip

\noindent{\it Proof of Theorem 1.1}. Let $m$ be the largest nonnegative integer not exceeding $n$
such that $\sum_{0<i\ls m}\lfloor(|A_i|-1)/k\rfloor<p(F)$.
For each $0<i\ls m$ let $A_i'$ be a subset of $A_i$ with cardinality $k\lfloor(|A_i|-1)/k\rfloor+1$.
In the case $m<n$, $p=p(F)$ is a prime and we let $A_{m+1}'$ be a subset of $A_{m+1}$
with
\begin{align*}|A_{m+1}'|=&k\(p-1-\sum_{0<i\ls m}\l\lfloor\f{|A_i|-1}k\r\rfloor\)+1
\\&<k\l\lfloor\f{|A_{m+1}|-1}k\r\rfloor+1\ls|A_{m+1}|.
\end{align*}
If $m+1<j\ls n$ then we let $A_{j}'\se A_j$ be a singleton.
Whether $m=n$ or not, we have $\sum_{i=1}^n(|A_i'|-1)=k(N-1),$
where
$$N=\min\l\{p(F),\,\sum_{i=1}^n\l\lfloor\f{|A_i|-1}k\r\rfloor+1\r\}.$$

Set
$$C=\{f(x_1,\ldots,x_n):\, x_1\in A_1',\ldots,x_n\in A_n'\}.$$
Suppose that $|C|\ls N-1$. Then
\begin{align*}&[x_1^{|A_1'|-1}\cdots x_n^{|A_n'|-1}]f(x_1,\ldots,x_n)^{N-1-|C|}\prod_{c\in C}(f(x_1,\ldots,x_n)-c)
\\=&[x_1^{|A_1'|-1}\cdots x_n^{|A_n'|-1}](a_1x_1^k+\cdots+a_nx_n^k)^{N-1}
\\=&\f{(N-1)!}{\prod_{i=1}^n((|A_i'|-1)/k)!}a_1^{(|A_1'|-1)/k}\cdots a_n^{(|A_n'|-1)/k}\not=0.
\end{align*}
By the Combinatorial Nullstellensatz, for some $x_1\in A_1',\ldots,x_n\in A_n'$
we have
$$f(x_1,\ldots,x_n)^{N-1-|C|}\prod_{c\in C}(f(x_1,\ldots,x_n)-c)\not=0$$
which contradicts the fact $f(x_1,\ldots,x_n)\in C$.

 In view of the above,
 $$|\{f(x_1,\ldots,x_n):\,x_1\in A_1,\ldots,x_n\in A_n\}|\gs|C|\gs N$$
 and this concludes the proof. \qed

\medskip
\noindent{\bf Example 2.1}. Let $a_1,\ldots,a_n\in F^*=F\sm\{0\}$, where $F$ is an algebraically closed field.
For each $i=1,\ldots,n$ let
$$A_i=\{x\in F:\, f_k(x)\in\{a_i^{-1},2a_i^{-1},\ldots,q_ia_i^{-1}\}\}\sm R_i,$$
where $k$ and $q_i<p(F)$ are positive integers,
\begin{equation} \label{eq-2.1}f_k(x)=
\begin{cases} x^k-x&\t{if}\ p(F)\ \t{is a prime divisor of}\ k,\\x^k&\t{otherwise},
\end{cases}
\end{equation}
and $R_i$ is a subset of $\{x\in F:\,f_k(x)=q_ia_i^{-1}\}$
with $|R_i|=r_i\ls k-1$.
For each $c\in F^*$, the equation $f_k(x)=c$ has exactly $k$ distinct solutions in $F$
since there is no $x\in F$ satisfying $f_k(x)-c=0=f'_k(x)$, where $f_k'$ is the formal derivative of $f_k$.
So $|A_i|=kq_i-r_i$ and hence
$$\l\lfloor\f{|A_i|-1}k\r\rfloor=\l\lfloor\f{k(q_i-1)+(k-1-r_i)}k\r\rfloor=q_i-1.$$
For every $i=1,\ldots,n$ we have
$$\{f_k(a):\,a\in A_i\}=\{a_i^{-1},2a_i^{-1},\ldots,q_ia_i^{-1}\}.$$
Thus
\begin{align*} &\{a_1f_k(x_1)+\cdots+a_nf_k(x_n):\, x_1\in A_1,\ldots,x_n\in A_n\}
\\=&\{a_1(y_1a_1^{-1})+\cdots+a_n(y_na_n^{-1}):\, y_i\in\{1,\ldots,q_i\}\ \t{for}\ i=1,\ldots,n\}
\\=&\{(y_1+\cdots+y_n)e:\, y_i\in\{1,\ldots,q_i\}\ \t{for}\ i=1,\ldots,n\}
\\=&\{ne,(n+1)e,\ldots,(q_1+\cdots+q_n)e\},
\end{align*}
where $e$ denotes the multiplicative identity of the field $F$.
Observe that
$$q_1+\cdots+q_n-n=\sum_{i=1}^n(q_i-1)=\sum_{i=1}^n\l\lfloor\f{|A_i|-1}k\r\rfloor.$$
Therefore
\begin{align*} &|\{a_1f_k(x_1)+\cdots+a_nf_k(x_n):\, x_1\in A_1,\ldots,x_n\in A_n\}|
\\&\quad=\min\bg\{p(F),\,\sum_{i=1}^n\l\lfloor\f{|A_i|-1}k\r\rfloor+1\bg\}.
\end{align*}

\section{Proofs of Theorem 1.2 and Corollary 1.4}\par

\noindent{\it Proof of Theorem 1.2}. Let $q_i=\lfloor(|A_i|-i)/k\rfloor$
for $i=1,\ldots,n$. And let $m$ be the largest nonnegative integer not exceeding $n$
such that $\sum_{0<i\ls m}q_i<p(F)$.
For each $0<i\ls m$ let $A_i'$ be a subset of $A_i$ with cardinality $kq_i+i$.
In the case $m<n$, $p=p(F)$ is a prime and we let $A_{m+1}'$ be a subset of $A_{m+1}$
with
\begin{align*}|A_{m+1}'|=&k\(p-1-\sum_{0<i\ls m}q_i\)+m+1
\\&<kq_{m+1}+m+1\ls|A_{m+1}|.
\end{align*}
If $m+1<j\ls n$ then we let $A_{j}'\se A_j$ with $|A_j|=j$. Whether $m=n$ or not, we have
$\sum_{i=1}^n(|A_i'|-i)=k(N-1),$ where
$$N=\min\{p(F),\,q_1+\cdots+q_n+1\}.$$

Set
$$C=\{f(x_1,\ldots,x_n):\, x_1\in A_1',\ldots,x_n\in A_n',\,\t{and}\ x_i\not=x_j\ \t{if}\ i\not=j\}.$$
Suppose that $|C|\ls N-1$ and let $P(x_1,\ldots,x_n)$ denote the polynomial
$$f(x_1,\ldots,x_n)^{N-1-|C|}\prod_{c\in C}(f(x_1,\ldots,x_n)-c)
\times\prod_{1\ls i<j\ls n}(x_j-x_i).$$
Then
$$\deg P\ls k(N-1)+\bi n2=\sum_{i=1}^n(|A_i'|-1).$$
By linear algebra,
$$\prod_{1\ls i<j\ls n}(x_j-x_i)=\det(x_i^{j-1})_{1\ls i,j\ls n}
=\sum_{\sigma\in S_n}\sign(\sigma)\prod_{i=1}^nx_i^{\sigma(i)-1},$$
where $S_n$ is the symmetric group of all permutations on $\{1,\ldots,n\}$,
and $\sign(\sigma)$ is $1$ or $-1$ according as $\sigma$ is even or odd.
Recall that
$$|A_i'|=
\begin{cases} kq_i+i&\t{if}\ 1\ls i\ls m,\\k(p(F)-1-q_1-\cdots-q_m)+i&\t{if}\ i=m+1\ls n,
\\i&\t{if}\ m+1<i\ls n.
\end{cases}$$
For each $i=1,\ldots,n$, we clearly have $|A_i'|-1\eq i-1\ (\mo\ k)$
and $\sigma(i)-1<n\ls k$ for all $\sigma\in S_n$.
Thus
\begin{align*}&\l[x_1^{|A_1'|-1}\cdots x_n^{|A_n'|-1}\r]P(x_1,\ldots,x_n)
\\=&\l[x_1^{|A_1'|-1}\cdots x_n^{|A_n'|-1}\r](a_1x_1^k+\cdots+a_nx_n^k)^{N-1}
\sum_{\sigma\in S_n}\sign(\sigma)\prod_{i=1}^nx_i^{\sigma(i)-1}
\\=&\bigg[\prod_{i=1}^nx_i^{|A_i'|-1}\bigg]
\sum_{j_1+\cdots+j_n=N-1}\f{(N-1)!}{j_1!\cdots j_n!}a_1^{j_1}\cdots a_n^{j_n}\prod_{i=1}^nx_i^{j_ik+i-1}
\\=&\f{(N-1)!}{q_1!\cdots q_m!q!}a_1^{q_1}\cdots a_m^{q_m}a\not=0,
\end{align*}
where $q=a=1$ if $m=n$, and $q=(p(F)-1-q_1-\cdots-q_m)!$ and $a=a_{m+1}^q$
if $m<n$.
In light of the Combinatorial Nullstellensatz, there are $x_1\in A_1',\ldots,x_n\in A_n'$
such that $P(x_1,\ldots,x_n)\not=0$. Obviously this contradicts the fact $f(x_1,\ldots,x_n)\in C$.

 By the above, $|V|\gs|C|\gs N$ and hence (1.5) holds. \qed

\medskip
\noindent{\bf Example 3.1}. Let $F$ be any algebraically closed field, and let $a_1,\ldots,a_n\in F^*$.
For each $i=1,\ldots,n$ let
$$A_i=S_i\cup\{x\in F:\, f_k(x)\in\{ja_i^{-1}:\,1<j\ls q_i\}\},$$
where $k\gs n$ and $q_i<p(F)$ are positive integers, $f_k(x)$ is given by (2.1),
and $S_i$ is a subset of $\{x\in F:\,f_k(x)=a_i^{-1}\}$
with $|S_i|\gs i$. Observe that $|A_i|=k(q_i-1)+|S_i|$ and hence
$$\l\lfloor\f{|A_i|-i}k\r\rfloor=q_i-1.$$
For every $i=1,\ldots,n$, we have
$$\{f_k(a):\,a\in A_i\}=\{a_i^{-1},2a_i^{-1},\ldots,q_ia_i^{-1}\}.$$
If $y_i\in\{1,\ldots,q_i\}$ for $i=1,\ldots,n$,
we can find distinct $x_1\in A_1,\ldots,x_n\in A_n$ such that $f_k(x_i)=y_ia_i^{-1}$ for $i=1,\ldots,n$;
in fact, if $x_1\in A_1,\ldots,x_{i-1}\in A_{i-1}$ are distinct with $i\ls n$, and
$f_k(x_j)=y_ja_j^{-1}$ for $j=1,\ldots,i-1$,
then we can choose $x_i\in A_i\sm\{x_1,\ldots,x_{i-1}\}$ satisfying $f_k(x_i)=y_ia_i^{-1}$
because $k\gs |S_i|>i-1$. Thus
\begin{align*} &\{a_1f_k(x_1)+\cdots+a_nf_k(x_n):\, x_1\in A_1,\ldots,x_n\in A_n,\ \t{and}\ x_i\not=x_j\ \t{if}\ i\not=j\}
\\=&\{a_1(y_1a_1^{-1})+\cdots+a_n(y_na_n^{-1}):\, y_i\in\{1,\ldots,q_i\}\ \t{for}\ i=1,\ldots,n\}
\\=&\{(y_1+\cdots+y_n)e:\, y_i\in\{1,\ldots,q_i\}\ \t{for}\ i=1,\ldots,n\}
\\=&\{ne,(n+1)e,\ldots,(q_1+\cdots+q_n)e\},
\end{align*}
where $e$ is the multiplicative identity of $F$.
Note that $$q_1+\cdots+q_n-n=\sum_{i=1}^n(q_i-1)=\sum_{i=1}^n\l\lfloor\f{|A_i|-i}k\r\rfloor.$$
So we have
\begin{align*} &|\{a_1f_k(x_1)+\cdots+a_nf_k(x_n):\, x_1\in A_1,\ldots,x_n\in A_n,
\ \t{and}\ x_i\not=x_j\ \t{if}\ i\not=j\}|
\\&\qquad\qquad=\min\bg\{p(F),\,\sum_{i=1}^n\l\lfloor\f{|A_i|-i}k\r\rfloor+1\bg\}.
\end{align*}

\bigskip

\noindent{\it Proof of Corollary 1.4}. The case $|A|<n$ is trivial; below we assume $|A|\gs n$.

We first handle the case $n\ls k$. If $1\ls i\ls\min\{n,r\}$ then
we let $A_i$ be a subset of $A$ with cardinality
$kq+i+\max\{r-n,0\}\ls kq+r=|A|$; if $r<j\ls n$ then we let $A_j$
be a subset of $A$ with cardinality $k(q-1)+j\ls k(q-1)+n\ls kq$.
(Note that when $r<n$ we have $q\not=0$ since $|A|\gs n$.)
Obviously,
$$\sum_{i=1}^n\l\lfloor\f{|A_i|-i}k\r\rfloor
=\sum_{i=1}^{\min\{n,r\}}q+\sum_{r<j\ls
n}(q-1)=n(q-1)+\min\{n,r\}.$$
Applying Theorem 1.2 we obtain that
\begin{align*} &|\{f(x_1,\ldots,x_n):\, x_1,\ldots,x_n\in A,\,\t{and}\ x_i\not=x_j\ \t{if}\ i\not=j\}|
\\\gs&|\{f(x_1,\ldots,x_n):\, x_1\in A_1,\ldots,x_n\in A_n,\,\t{and}\ x_i\not=x_j\ \t{if}\ i\not=j\}|
\\\gs&\min\{p(F),\,n(q-1)+\min\{n,r\}+1\}.
\end{align*}
In particular, if $k=n$ then
\begin{align*} &|\{f(x_1,\ldots,x_n):\, x_1,\ldots,x_n\in A,\,\t{and}\ x_i\not=x_j\ \t{if}\ i\not=j\}|
\\\gs&\min\{p(F),\,n(q-1)+r+1\}=\min\{p(F),\,|A|-n+1\}.
\end{align*}

Now suppose that $n>k$. Let $c_{k+1},\ldots,c_n$ be $n-k$ distinct elements of $A$.
Then $A'=A\sm\{c_{k+1},\ldots,c_n\}$ has cardinality $|A|-(n-k)$. By what we have proved,
\begin{align*} &|\{f(x_1,\ldots,x_k,c_{k+1},\ldots,c_n):\, x_1,\ldots,x_k\in A',
\,\t{and}\ x_i\not=x_j\ \t{if}\ i<j\}|
\\&\qquad\gs\min\{p(F),\,|A'|-k+1\}=\min\{p(F),\,|A|-n+1\}.
\end{align*}
So the desired inequality follows. \qed

\section*{4. Discussion of Conjecture 1.1}\par

Conjecture 1.1 in the case $f(x_1,\ldots,x_n)=x_1+\cdots+x_n$, essentially gives
the Dias da Silva--Hamidoune result mentioned in Section 1.

In the case $f(x_1,\ldots,x_n)=x_1^k+\cdots+x_n^k$ with $k>1$, we may explain
the symmetry between $n$ and $|A|-n$ as follows: If $|A|>n$ then
\begin{align*}&|\{x_1^k+\cdots+x_n^k:\, x_1,\ldots,x_n\in A,\,\t{and}\ x_i\not=x_j\ \t{if}\ i\not=j\}|
\\=&\bg|\bg\{\sum_{a\in A}a^k-y_1^k-\cdots-y_{|A|-n}^k:
\, y_1,\ldots,y_{|A|-n}\in A,\,\t{and}\ y_i\not=y_j\ \t{if}\ i\not=j\bg\}\bg|
\\=&|\{y_1^k+\cdots+y_{|A|-n}^k:\, y_1,\ldots,y_{|A|-n}\in A,\,\t{and}\ y_i\not=y_j\ \t{if}\ i\not=j\}|.
\end{align*}

Conjecture 1.1 holds when $n=2$. In fact, if we set $A_1=\{a_1x:\,x\in A\}$ and $A_2=\{a_2x:\,x\in A\}$,
then $|A_1|+|A_2|-3=2(|A|-2)+1$ and
\begin{align*}&|\{a_1x_1+a_2x_2:\, x_1,x_2\in A\ \t{and}\ x_1\not=x_2\}|
\\=&|\{y_1+y_2:\, y_1\in A_1,\ y_2\in A_2\ \t{and}\ a_1^{-1}y_1\not=a_2^{-1}y_2\}|
\\=&|\{y_1+y_2:\, y_1\in A_1,\ y_2\in A_2\ \t{and}\ y_1-a_1a_2^{-1}y_2\not=0\}|
\\\gs&\min\{p(F)-\da,\,|A_1|+|A_2|-3\}\ \ (\t{by Corollary 3 of \cite{PS}}).
\end{align*}

The following example illustrates how the lower bound in (1.7) comes out.
\medskip

\noindent{\bf Example 4.1}. Let $k,n\in\Z^+$ and $q\in\N$.
Let $$A=\{z\in\C:\, z^k\in\{1,\ldots,q\}\}\cup R,$$
where $\C$ is the field of complex numbers and $R$ is a subset of $\{z\in\C:\,z^k=q+1\}$
with cardinality $r<k$. Suppose that $n\gs k$ and
$|A|=kq+r\gs n=k\lfloor n/k\rfloor+s$ with $s\in\{0,\ldots,k-1\}$.
Clearly,
$$V=\{x_1^k+\cdots+x_n^k:\, x_1,\ldots,x_n\in A,\, \t{and}\ x_i\not=x_j\ \t{if}\ i\not=j\}$$
just consists of those sums $\sum_{i=1}^ny_i$ with
$y_i\in\{1,\ldots,q+1\}$, $|\{i:\,y_i=q+1\}|\ls r$ and
$|\{i:\,y_i=j\}|\ls k$ for all $j=1,\ldots,q$. Thus the smallest element of $V$ is
$$m_V=k\sum_{i=1}^{\lfloor n/k\rfloor}i+s\l(\l\lfloor\f nk\r\rfloor+1\r)
=\l(\f k2\l\lfloor\f nk\r\rfloor+s\r)\l(\l\lfloor\f nk\r\rfloor+1\r),$$
while the largest element of $V$ is
$$M_V=r(q+1)+k\sum_{i=0}^{\lfloor n/k\rfloor-1}(q-i)+(s-r)\l(q-\l\lfloor \f nk\r\rfloor\r)-d(r,s),$$
where
$$d(r,s)=
\begin{cases} r-s&\t{if}\ r\gs s,\\0&\t{if}\ r<s.
\end{cases}$$
It follows that
\begin{align*} &d(r,s)+M_V-m_V
\\=&r(q+1)+kq\l\lfloor\f nk\r\rfloor-k\l\lfloor \f nk\r\rfloor^2
+(s-r)q+(r-s)\l\lfloor\f nk\r\rfloor-s\l(\l\lfloor\f nk\r\rfloor+1\r)
\\=&\l\lfloor\f nk\r\rfloor\l(kq+r-s-k\l\lfloor \f nk\r\rfloor-s\r)+r+sq-s
\\=&\l\lfloor\f nk\r\rfloor(|A|-s-n)+s(q-1)+r
\\=&\l\lfloor\f nk\r\rfloor(|A|-n)+s\l(q-\l\lfloor\f nk\r\rfloor-1\r)+r
\\=&\l\lfloor\f nk\r\rfloor(|A|-n)+s\l\lfloor\f{|A|-n}k\r\rfloor+
\begin{cases} r-s&\t{if}\ r\gs s,\\r&\t{if}\ r<s.
\end{cases}
\end{align*}
Therefore
\begin{align*} |V|=&|\{m_V,\ldots,M_V\}|=M_V-m_V+1
\\=&1+\l\lfloor\f nk\r\rfloor(|A|-n)+s\l\lfloor\f{|A|-n}k\r\rfloor+
\begin{cases} r&\t{if}\ r<s,\\0&\t{otherwise}.
\end{cases}
\end{align*}
Finally, we note that
\begin{align*}&\l\lfloor\f nk\r\rfloor(|A|-n)+s\l\lfloor\f{|A|-n}k\r\rfloor
\\=&\l(\f nk-\l\{\f nk\r\}\r)(|A|-n)+k\l\{\f nk\r\}\l\lfloor\f{|A|-n}k\r\rfloor
\\=&\f{n(|A|-n)}k-k\l\{\f nk\r\}\l\{\f{|A|-n}k\r\}.
\end{align*}

\Ack. The author thanks the referees for their helpful comments.

\end{document}